\documentclass[a4paper,twoside,12pt]{article}
\usepackage{graphicx,psfig,epsfig}
\usepackage{amsfonts}
\usepackage{amssymb}
\usepackage{amsmath}
\usepackage{subfigure}
\usepackage{setspace}
\usepackage[text={16cm,24cm},top=3cm,bottom=3cm,centering]{geometry}
\renewcommand{\maketitle}%
{\thispagestyle{empty}
\noindent
\par\vspace*{3cm}\par
{\centering{\Large\textbf\MATCHtitle\par}
\par\vspace*{12pt}\par
{\textrm\MATCHauthor}\par\vspace*{12pt}\par}%
}
\renewenvironment{abstract}%
{\begin{center}\begin{minipage}[t]{1.0\textwidth}
\noindent\textbf{Abstract}$\colon$}%
{\end{minipage}\end{center}\par}	 

\usepackage{enumerate}
\usepackage[square,sort]{natbib}
\usepackage[nosepfour,warning,np,debug,autolanguage]{numprint}
\newenvironment{classification}%
{\begin{center}\begin{minipage}[t]{1.0\textwidth}
\noindent\textbf{AMS Mathematics Subject Classification (2000)}$\colon$}%
{\end{minipage}\end{center}}
\newenvironment{keywords}%
{\begin{center}\begin{minipage}[t]{1.0\textwidth}
}%
{\end{minipage}\end{center}\par}
\newcommand{\MATCHtitle}%
{Grid Generation and Adaptation by Functionals}
\newcommand{\MATCHauthor}
{\text{Sanjay Kumar Khattri}
\\[1ex]
\text{Department of Mathematics, The University of Bergen, Norway}\\
\tt{sanjay@mi.uib.no}\\ \texttt{{http://www.mi.uib.no/$\sim$sanjay}}}
\begin{document}
\begin{spacing}{1.0}
\maketitle
\begin{abstract}
Accuracy of a simulation is strongly depend on the grid quality. Here, quality means
orthogonality at the boundaries and quasi-orthogonality within the critical regions,
smoothness, bounded aspect ratios, solution adaptive behaviour, etc. We review various functionals for generating high quality structured quadrilateral meshes in two dimensional domains. Analysis of Winslow and Modified Liao functionals are presented. Numerical examples are also presented to support our theoretical analysis. We will demonstrate the use of the Area functional for generating adaptive quadrilateral meshes. 
\end{abstract}
\section{Introduction}
\label{sec:intro}
Accuracy of numerical solutions of partial differential equations on a grid is very much depend on the quality of the underlying grid. There are various parameters for measuring grid quality. For example, orthogonality of grid lines and grid density in the regions of large solution gradients. A desired grid may be an orthogonal grid with high grid density in the areas of sharp solution gradients. Variational methods has been used for improving quality of a given grid \citep{gridbook}. In the variational methods, a grid functional is def{i}ned.  Grid functional is an algebraic expression of the position vectors of the internal nodes of a mesh. Optimization of the grid functional may result in a grid with desired properties such as orthogonal grid lines, equal cell areas, linear or parallelogram cells \citep[see][]{khattri_1} and untangled mesh \citep{knupp1,knupp2,knupp3}. There are many algebraic functionals for grid generation and optimization \citep[cf.][]{Tinoco1,Tinoco2,knupp1,knupp2,khattri_1,knupp3,Tinoco3}. The first study of grid generation by algebraic functionals were done in \citep{Castillo}. Castillo and Steinberg introduced Length, Orthogonality and Area functionals \citep{Castillo}. Area functional are well known for producing robust quadrilateral meshes. For a detailed description of properties of area functionals, the interested readers are referred to \citep{Tinoco3}. Recently the area functional has been used for generating adaptive quadrilateral meshes \citep{khattri_2}. 

Let $x(\xi,\eta)$ and $y(\xi,\eta)$ be the coordinates of a node in a mesh. Let us further assume that $x$ and $y$ are twice differentiable functions of the independent variables $\xi$ and $\eta$. An integral functional $\mathcal{I}$ can be defined as follows
\begin{equation}
\mathcal{I}(x,y) \overset{\textbf{def}}{=} \int_{[0,1]\times[0,1]}{\mathcal{F}(\xi,\eta,x,y,x_\xi,x_\eta,y_\xi,y_\eta)}\,d\xi\,d\eta\enspace{.}
\end{equation}
We are interested in finding the functions $x(\xi,\eta)$ and $y(\xi,\eta)$ for which the integral functional $\mathcal{I}$ attains an extremal value. Such coordinates $x$ and $y$ define a mesh with desirable properties. The integral functional $\mathcal{I}$ is also referred to as control function for adaptive grid generation \citep{gridbook}. The conditions for the extremal value of the integral functional $\mathcal{I}$ are expressed by the Euler-Lagrange equations. The two Euler-Lagrange equations are
\begin{alignat}{2}
\dfrac{\partial\mathcal{F}}{\partial{x}}-\dfrac{\partial}{\partial{\xi}}\left({\dfrac{\partial\mathcal{F}}{\partial{x_\xi}}}\right)-\dfrac{\partial}{\partial{\eta}}\left({\dfrac{\partial\mathcal{F}}{\partial{x_\eta}}}\right) = 0\enspace{,} \\
\dfrac{\partial\mathcal{F}}{\partial{y}}-\dfrac{\partial}{\partial{\xi}}\left({\dfrac{\partial\mathcal{F}}{\partial{y_\xi}}}\right)-\dfrac{\partial}{\partial{\eta}}\left({\dfrac{\partial\mathcal{F}}{\partial{y_\eta}}}\right) = 0\enspace{.}
\end{alignat}
\begin{figure}
\begin{minipage}[t]{6.3cm}
\begin{center}
\includegraphics[width=1.3\textwidth]{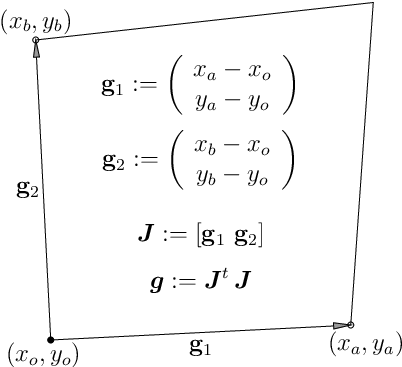}
\caption{Quantities of interest for a quadrilateral cell.}
\label{fig:quad_cell}
\end{center}
\end{minipage} 
\hfill
\begin{minipage}[t]{6.3cm}
\begin{center}
\includegraphics[width=1.3\textwidth]{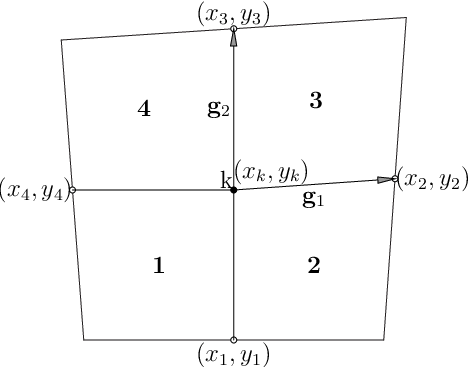}
\caption{{2D Structured Mesh. Node $k$ is surrounded by four quadrilaterals.}}
\label{fig:2d_mesh}
\end{center}
\end{minipage}
\hfill
\end{figure}
The functions $x$ and $y$, which satisfy the above Euler-Lagrangian equations, are called the extremals of the integral functional $\mathcal{I}$. 
 
Let us def{i}ne some quantities of interest. Figure \ref{fig:quad_cell} shows a quadrilateral cell, and this cell belongs to a mesh. Let the co-variant vector at the node $o$ and in the direction $oa$ is $\mathbf{g}_1$, and another co-variant vector at the node $o$ but in the direction $ob$ is $\mathbf{g}_2$. These vectors are given as
\begin{equation}
\mathbf{g}_1 = (x_a-x_o,y_a-y_o)^t \quad \text{and}\quad  \mathbf{g}_2 = (x_b-x_o,y_b-y_o)^t\enspace{.}
\end{equation}
Other interesting quantities such as the Jacobian and g-tensor matrix can be def{i}ned from the co-variant vectors. The columns of the Jacobian matrix are the co-variant vectors. The g-tensor matrix is the product of the Jacobian matrix with it's transpose. Thus, the Jacobian matrix and the g-tensor at the node $o$ and for the cell shown in the Figure \ref{fig:quad_cell} are given as
\begin{equation}
\boldsymbol{J} = \left[\mathbf{g}_1\,\,\mathbf{g}_2\right] \quad \text{and}\quad \boldsymbol{g} = \boldsymbol{J}^t\,\boldsymbol{J}\enspace{.}
\end{equation}

The layout of the paper is as follows. In the Section \ref{discrete_fun}, several functionals are presented. Continuous and discrete versions of the functionals are presented. Section \ref{num_examp} presents several numerical experiments, and finally Section \ref{conclusion_} concludes the paper.
\section{Discrete Functionals}
\label{discrete_fun}
Let us f{i}rst introduce some quantities of interest. These will be used later in formulating algebraic functionals. 
Figure \ref{fig:2d_mesh} is a $2\times2$ structured mesh. We use this f{i}gure for def{i}ning these quantities. 

$J(k_i)$ refers to the Jacobian (determinant of the Jacobian matrix $\boldsymbol{J}(k_i)$) at the node $k$ and for the cell $i$. Table \ref{jacobian_table} lists the Jacobian matrix for the four cells surrounding the node $k$. $\mathbf{g}_1(k_i)$ refers to the co-variant base vector at the node $k$ and for the cell $i$. The base vector $\mathbf{g}_1$  points along horizontal grid lines. Similarly, $\mathbf{g}_2(k_i)$ refers to the co-variant base vector at the node $k$ and for the cell $i$, and it points along the vertical grid lines. Table \ref{covec_table} lists the co-variant vectors for the Figure \ref{fig:2d_mesh}. It should be noted that column vectors of the Jacobian matrix are the co-variant base vectors. For example, the column vectors of $\boldsymbol{J}(k_1)$ are $\mathbf{g_1}{(k_1)}$ and $\mathbf{g_2}{(k_1)}$. That is $\boldsymbol{J}(k_1)=\left[\mathbf{g}_1{(k_1)}\,\,\mathbf{g}_2{(k_1)}\right]$.

$\boldsymbol{g}(k_i)$ refers to the co-variant metric tensor at the node $k$ and for the cell $i$. It is def{i}ned as $\boldsymbol{g}(k_i) = {\boldsymbol{J}(k_i)}^t\,\boldsymbol{J}(k_i)$. $g_{mn}(k_i)$ refers to the $(m, n)$ coeff{i}cient of the co-variant metric tensor $\boldsymbol{g}(k_i)$ for the node $k$ and for the cell $i$. It can be seen that $g_{11}(k_i)={\mathbf{g}_{1}(k_i)}^t\cdot\mathbf{g}_{1}(k_i)$ and $g_{12}(k_i)={\mathbf{g}_{1}(k_i)}^t\cdot\mathbf{g}_{2}(k_i)$. Similarly, other coeff{i}cients can be def{i}ned. 

The coeff{i}cient $g_{ 12}$ is a measure of the angle between the co-variant base vectors $\mathbf{g}_1$ and $\mathbf{g}_2$. While, the coeff{i}cient $\mathbf{g}_{11}$ is a measure of the discrete $L_2$ length of the co-variant vector $\mathbf{g}_1$.
\begin{table}
\hfill
\caption{Jacobian matrix at the node $k$ for the surrounding cells for the Figure \ref{fig:2d_mesh}.}
\centering
\begin{tabular}{ccccc}
\hline\\
$\boldsymbol{J}{(k_1)}=\left[
\begin{array}{cc}
(x_4-x_k) & (x_1-x_k) \\
(y_4-y_k) & (y_1-y_k) 
\end{array}
 \right]$& $\boldsymbol{J}(k_2)=\left[
\begin{array}{cc}
(x_2-x_k) & (x_1-x_k) \\
(y_2-y_k) & (y_1-y_k) 
\end{array}
 \right]$ \\ \\ 
  $\boldsymbol{J}(k_3)=\left[
\begin{array}{cc}
(x_2-x_k) & (x_3-x_k) \\
(y_2-y_k) & (y_3-y_k) 
\end{array}
 \right]$ &
 $\boldsymbol{J}(k_4)=\left[
\begin{array}{cc}
(x_4-x_k) & (x_3-x_k) \\
(y_4-y_k) & (y_3-y_k) 
\end{array}
 \right]$  \\ `
  \\
\hline \hline
\end{tabular}
\label{jacobian_table}
\end{table}
\begin{table}
\hfill
\caption{Co-variant vectors at the node $k$ for the surrounding cells for the Figure \ref{fig:2d_mesh}.}
\begin{tabular}{ccccc}
\hline\\
$\mathbf{g}_1(k_1)=
\left(
\begin{array}{c}
x_4-x_k \\ 
y_4-y_k
\end{array}
\right)$ & 
$\mathbf{g}_2(k_1)=
\left(
\begin{array}{c}
x_1-x_k \\ 
y_1-y_k
\end{array}
\right)$ & 
$\mathbf{g}_1(k_2)=
\left(
\begin{array}{c}
x_2-x_k \\ 
y_2-y_k
\end{array}
\right)$  & 
$\mathbf{g}_2(k_2)=
\left(
\begin{array}{c}
x_1-x_k \\ 
y_1-y_k
\end{array}
\right)$ \\ \\ 
$\mathbf{g}_1(k_3)=
\left(
\begin{array}{c}
x_2-x_k \\ 
y_2-y_k
\end{array}
\right)$ & 
$\mathbf{g}_2(k_3)=
\left(
\begin{array}{c}
x_3-x_k \\ 
y_3-y_k
\end{array}
\right)$ & 
$\mathbf{g}_1(k_4)=
\left(
\begin{array}{c}
x_4-x_k \\ 
y_4-y_k
\end{array}
\right)$ & 
$\mathbf{g}_2(k_4)=
\left(
\begin{array}{c}
x_3-x_k \\ 
y_3-y_k
\end{array}
\right)$ \\ \\
\hline \hline
\end{tabular}
\label{covec_table}
\end{table}

Let us consider a structured quadrilateral mesh (each internal node is surrounded by four quadrilaterals) consisting of $n$ internal nodes. The following functionals can be defined
\subsection{Area Functional}
The integral form of the standard Area functional is given as
\begin{alignat}{2}
\mathcal{I}_\text{A} &\overset{\textbf{def}}{=} \dfrac{1}{2} \,\int_{[0,1]\times[0,1]}{\vert{\boldsymbol{J}}\vert}^2\,d\xi\,d\eta\enspace{,} \\
&=\int_{[0,1]\times[0,1]}({x_\xi\,y_\eta-x_\eta\,y_\xi})\,d\xi\,d\eta\enspace{.}
\end{alignat}
The Euler-Lagrangian equations for the Area functional are
\begin{alignat}{2}
\dfrac{\partial}{\partial{\xi}}(\vert{\boldsymbol{J}}\vert\,x_\eta) - \dfrac{\partial}{\partial{\eta}}(\vert{\boldsymbol{J}}\vert\,x_\xi) &= 0\enspace{,}\\
\dfrac{\partial}{\partial{\xi}}(\vert{\boldsymbol{J}}\vert\,y_\eta) - \dfrac{\partial}{\partial{\eta}}(\vert{\boldsymbol{J}}\vert\,y_\xi) &= 0\enspace{.}
\end{alignat}
In the simplif{i}ed form the above equations can be written as
\begin{alignat}{2}
{y_\eta}^2\,{x_{\xi\xi}}-{x_\eta}\,y_\eta\,{y_{\xi\xi}}-2.0\,y_{\xi}\,y_\eta\,x_{\xi\eta}+(x_\xi\,y_\eta+x_\eta\,y_\xi)\,y_{\xi\eta}+{y_\xi}^2\,x_{\eta\eta}-x_\xi\,y_\xi\,y_{\eta\eta} = 0\enspace{,}\\
{x_\eta}^2\,{y_{\xi\xi}}-{x_\eta}\,y_\eta\,{x_{\xi\xi}}-2.0\,x_{\xi}\,x_\eta\,y_{\xi\eta}+(x_\xi\,y_\eta+x_\eta\,y_\xi)\,x_{\xi\eta}+{x_\xi}^2\,y_{\eta\eta}-y_\xi\,y_\xi\,x_{\eta\eta} = 0\enspace{,}
\end{alignat}
\citep[see][]{Tinoco3}.
The above Euler-Lagrangian equations are non-elliptic, coupled and quasi-linear \citep[cf.][]{Tinoco3}. 
For generating adaptive mesh, the author proposed the following variation in the Area functional  
\begin{equation}
\mathcal{F}_A (\mathbf{\mathbf{x},\mathbf{y}}) = \sum_{k=1}^n\left[{\sum_{i=1}^4{s(k_i)\left[J(k_i)\right]^2}}\right]\enspace{,}
\label{area_fun_1}
\end{equation}
\citep{khattri_2}. In the above equation, $s(k)$ is called the adaptive function, and $s(k_i)$ is the value of the adaptive function at the node $k$ and for cell $i$.
\subsection{Length Functional}
The integral form of the Length functional is given as 
\begin{alignat}{2}
\mathcal{I}_\text{L} &\overset{\textbf{def}}{=} \dfrac{1}{2} \int_{[0,1]\times[0,1]}\left[\,{g_{11}+g_{22}}\,\right]\,d\xi\,d\eta\enspace{,} \\
&=\dfrac{1}{2}\,{\int_{[0,1]\times[0,1]}{ \left[{(x_\xi)^2+{(x_\eta)^2}+(y_\xi)^2+(y_\eta)^2}\right]\,d\xi\,d\eta } }\enspace{,}
\end{alignat}
\citep[and references therein]{Tinoco1,Tinoco2,Tinoco3}. The conditions of extremality of the above length functional are given by the following Euler-Lagrangian equations
\begin{alignat}{2}
\dfrac{\partial^2{x}}{\partial{\xi}^2}+\dfrac{\partial^2{x}}{\partial{\eta}^2} = 0\enspace{,}\\
\dfrac{\partial^2{y}}{\partial{\xi}^2}+\dfrac{\partial^2{y}}{\partial{\eta}^2} = 0\enspace{.}
\end{alignat}
The above Laplace's equations can be solved in the computational domain $[0,1]\times[0,1]$ with a specified value of $x$ and $y$ on the boundary. The Euler-Lagrangian equations associated with the Length functional are linear and uncoupled.

The discrete Length functional \cite{Castillo,gridbook} is give as follows
\begin{equation}
\mathcal{F}_\text{L}(\mathbf{\mathbf{x},\mathbf{y}}) = \sum_{k=1}^n\left[{\sum_{i=1}^4\left({{g_{11}(k_i)}+{g_{22}(k_i)}}\right)
}\right]\enspace{.}
\label{area_fun}
\end{equation}
\subsection{Orthogonality Functional}
The integral form of the Orthogonality functional \citep[and references therein]{Tinoco1,Tinoco2,Tinoco3} is given as follows
\begin{alignat}{3}
\mathcal{I}_\text{O} &\overset{\textbf{def}}{=}  \dfrac{1}{2}\,\int_{[0,1]\times[0,1]}{(g_{12})^2}\,d\xi\,d\eta\enspace{,}\\ &=\dfrac{1}{2}\,\int_{[0,1]\times[0,1]}{(\mathbf{g}_1\cdot\mathbf{g}_2)^2}\,d\xi\,d\eta, \\ &=\dfrac{1}{2}\,\int_{[0,1]\times[0,1]} {(x_\xi\,x_\eta+y_\xi\,y_\eta)^2}\,d\xi\,d\eta\enspace{.}
\end{alignat}
The Euler-Lagrangian equations corresponding to the minimization of the above integral are 
\begin{alignat}{2}
\dfrac{\partial}{\partial\xi}\left(g_{12}\dfrac{\partial{x}}{\partial\eta}\right) + \dfrac{\partial}{\partial\eta}\left(g_{12}\dfrac{\partial{x}}{\partial\xi}\right) =0\enspace{,}\\
\dfrac{\partial}{\partial\xi}\left(g_{12}\dfrac{\partial{y}}{\partial\eta}\right) + \dfrac{\partial}{\partial\eta}\left(g_{12}\dfrac{\partial{y}}{\partial\xi}\right) =0\enspace{.}
\end{alignat}
These Euler-Lagrangian equations are quasilinear, coupled and non-elliptic in nature \cite{Tinoco3}. A simplif{i}ed form the above Euler-Lagrangian equations is
\begin{alignat}{2}
{x_\eta}^2\,x_{\xi\xi}+x_{\eta}\,y_{\eta}\,y_{\xi\xi}+(4\,x_\xi\,x_\eta+2\,y_\xi\,y_\eta)\,x_{\xi\eta}+(x_\xi\,y_\eta+x_\eta\,y_\xi)\,y_{\xi\eta}+{x_\xi}^2\,x_{\eta\eta}+x_{\xi}\,y_{\xi}\,y_{\eta\eta}=0\enspace{,}\\
{y_\eta}^2\,y_{\xi\xi}+y_{\eta}\,x_{\eta}\,x_{\xi\xi}+(4\,y_\xi\,y_\eta+2\,x_\xi\,x_\eta)\,y_{\xi\eta}+(y_\xi\,x_\eta+y_\eta\,x_\xi)\,x_{\xi\eta}+{y_\xi}^2\,y_{\eta\eta}+y_{\xi}\,x_{\xi}\,x_{\eta\eta}=0\enspace{,}
\end{alignat}
\citep[see][]{Tinoco3}. This functional takes only non-negative values, and it would attain a minimum value of zero for a completely orthogonal grid. The discrete version of the above Orthogonality functional \citep{Castillo,gridbook} is given as follows
\begin{equation}
\mathcal{F}_\text{O}(\mathbf{\mathbf{x},\mathbf{y}}) = \sum_{k=1}^n\left[{\sum_{i=1}^4{(\mathbf{g_1(k_i)}}\cdot{\mathbf{g_2(k_i)}})^2}\right]\enspace{.}
\label{ortho-1}
\end{equation}
It is found \citep[cf.][]{Tinoco1,Tinoco2,gridbook,Castillo,combi_1,combi_2} that a linear combination of Area, Length and Orthogonality functionals can produce robust grids in complicated 2D domains. 
\subsection{Combination of Length, Area and Orthogonality Functionals}
A combined functional is given as
\begin{equation}
\mathcal{F}(\mathbf{x},\mathbf{y}) = \Bbbk_\text{A}\,\mathcal{F}_\text{A}(\mathbf{x},\mathbf{y})+\Bbbk_\text{L}\,\mathcal{F}_\text{L}(\mathbf{x},\mathbf{y}) + \Bbbk_\text{O}\,\mathcal{F}_\text{O}(\mathbf{x},\mathbf{y})\enspace{,}
\end{equation}
\cite{Tinoco1,Tinoco2,gridbook,Castillo,combi_1,combi_2}. Here, the parameters $\Bbbk_\text{A}$, $\Bbbk_\text{L}$ and $\Bbbk_\text{O}$ satisfy$\colon$ $\Bbbk_\text{A}+\Bbbk_\text{L}+\Bbbk_\text{O} = 1.0$ and $\Bbbk_\text{A}\ge0$, $\Bbbk_\text{L}\ge0$, $\Bbbk_\text{O}\ge0$. A serious drawback of the above combined functional is a suitable choice of the parameters. It requires an experience in coming up with a good set of parameters \citep{Castillo}. It was found \cite{combi_1,combi_2} that the following choice of parameters
\begin{equation}
\Bbbk_\text{A} =0.50, \quad \Bbbk_\text{L} = 0.0, \quad \text{and}\quad \Bbbk_\text{O} = 0.50,
\end{equation}
produces robust grid in many practical domains. The corresponding functional is referred as the Knupp's functional. Presented numerical work shows that this functional can produce good grids. The Euler-Lagrangian \citep{Castillo} equations for the minimization of the Knupp's functional are
\begin{alignat}{2}
({x_\eta}^2+{y_\eta}^2)\,x_{\xi\xi} + 4\,x_\xi\,x_\eta\,x_{\xi\eta}+2(x_\xi\,y_\eta+x_\eta\,y_\xi)\,y_{\xi\eta}+({x_\xi}^2+{y_\xi}^2)\,x_{\eta\eta} &=0\enspace{,}\\
({x_\eta}^2+{y_\eta}^2)\,y_{\xi\xi} + 4\,y_\xi\,y_\eta\,y_{\xi\eta}+2(x_\xi\,y_\eta+x_\eta\,y_\xi)\,x_{\xi\eta}+({x_\xi}^2+{y_\xi}^2)\,y_{\eta\eta} &=0\enspace{.}
\end{alignat}
\subsection{Winslow Functional}
The Winslow functional is given as follows
\begin{equation}
\mathcal{F}(\mathbf{\mathbf{x},\mathbf{y}}) = \sum_{k=1}^n\left[{\sum_{i=1}^4\left({\dfrac{g_{11}(k_i)+g_{22}(k_i)}{{\vert \boldsymbol{J}(k_i)\vert}}}\right)}\right]\enspace{,}
\label{winslow}
\end{equation}
\citep{winslow,gridbook,shashkov1}. Here, ${\vert \boldsymbol{J}(k_i)\vert}$ is the determinant of the Jacobian matrix. One very important feature of the above functional is that it has barrier. It means the value of the functional approaches inf{i}nity when the cells degenerate. That is $\vert{\boldsymbol{J}}\vert\,\rightarrow\,0$. Thus, this functional produces unfolded grids. Numerical experiments also prove this feature of the Winlow functional. Since $g_{11}=\mathbf{g}_1\cdot \mathbf{g}_{1}$, and $g_{22}=\mathbf{g}_2\cdot \mathbf{g}_{2}$. It can be shown that the numerator ($g_{11}+g_{22}$) in the Winslow functional \eqref{winslow} is the Frobenius norm of the Jacobian matrix. That is $$g_{11}(k_i)+g_{22}(k_i) = \sum_{n=1}^2{\sum_{m=1}^2 (J_{mn}(k_i)})^2 = (\Vert{\mathbf{J(k_i)}}\Vert)^2\enspace{,}$$ 
Here, $J_{mn}$ are the components of the Jacobian matrix $\boldsymbol{J}$. Thus, the Winslow functional \eqref{winslow} can be written as follows
\begin{equation}
\mathcal{F}(\mathbf{\mathbf{x},\mathbf{y}}) = \sum_{k=1}^n\left[{\sum_{i=1}^4{\dfrac{\Vert{\boldsymbol{J}(k_i)}\Vert^2}{{\vert \boldsymbol{J}(k_i)\vert}}}}\right]\enspace{.}
\label{winslow2}
\end{equation}
It can be seen easily that the Frobenius norm a $2\times2$ matrix $\boldsymbol{A}$, and its inverse are related as $\Vert{\boldsymbol{A}^{-1}}\Vert=\dfrac{\Vert{\boldsymbol{A}}\Vert}{\vert{\boldsymbol{A}}\vert}$. Here, ${\vert{\boldsymbol{A}}\vert}$ is the determinant of the matrix $\boldsymbol{A}$. The condition number $\mathsf{K}(\boldsymbol{A})$ of a matrix $\boldsymbol{A}$ can be written as $\mathsf{K}(\boldsymbol{A})=\Vert{\boldsymbol{A}}\Vert\,\Vert{\boldsymbol{A}^{-1}}\Vert$. Here, the norm is the Frobenius norm. Thus, the Winslow functional can be written as follows
\begin{equation}
\mathcal{F}(\mathbf{\mathbf{x},\mathbf{y}}) = \sum_{k=1}^n\left[{\sum_{i=1}^4{\mathsf{K}(\boldsymbol{J}(k_i))}}\right]\enspace{.}
\label{winslow3}
\end{equation}
Thus, the minimization of the functional \eqref{winslow} is equivalent to the minimization of the condition number of the Jacobian matrix. A detailed description of the above analysis can also be found in \citep{knupp1,knupp2,knupp3,shashkov1}. The condition number $\mathsf{K}(\boldsymbol{J}(k_i))$ can also be expressed as $$\mathsf{K}(\boldsymbol{J}(k_i))=\dfrac{\mathbf{g_1}{(k_i)}^2+\mathbf{g_2}{(k_i)}^2}{\vert{\mathbf{g_1}{(k_i)}\times\mathbf{g_2}{(k_i)}}\vert}\enspace{.}$$
The $\boldsymbol{g}(k_i)$ tensor matrix is give as 
$$ 
\boldsymbol{g}(k_i)=\left[
\begin{array}{cc}
g_{11}(k_i) & g_{12}(k_i) \\ 
g_{21}(k_i) & g_{22}(k_i) 
\end{array}
\right]\enspace{.}$$
Let $\lambda_1$ and $\lambda_2$ be the eigenvalues of the matrix $\boldsymbol{g}(k_i)$. Then 
$${\dfrac{g_{11}(k_i)+g_{22}(k_i)}{{\vert \boldsymbol{J}(k_i)\vert}}}=\dfrac{\lambda_1+\lambda_2}{\sqrt{\lambda_1\,\lambda_2}}\ge2.0\enspace{.}$$
Here, we have used the relation $\vert{\boldsymbol{J}}\vert^2=\vert\boldsymbol{g}\vert$. Thus, the Winslow functional is bounded from below.
\subsection{Liao Functional}
The Liao functional for grid generation was proposed in \citep{liao}, and is give as follows
\begin{equation}
\mathcal{F}(\mathbf{\mathbf{x},\mathbf{y}}) = \sum_{k=1}^n\left[{\sum_{i=1}^4({{g_{11}}^2+{g_{22}}^2+2\,{g_{12}}^2})}\right]\enspace{.}
\label{liao}
\end{equation}
\subsection{Modified Liao Functional}
The Liao functional can produce folded grids. We will explore it through numerical experiments. In the literature, following modification \citep{gridbook} of the Liao functional is given
\begin{equation}
\mathcal{F}(\mathbf{\mathbf{x},\mathbf{y}}) = \sum_{k=1}^n\left[{\sum_{i=1}^4\left(\dfrac{g_{11}(k_i)+g_{22}(k_i)}{\sqrt{g(k_i)}}\right)^2}\right]\enspace{.}
\label{modliao}
\end{equation}
In the above equation, $g(k_i)$ is the determinant of the covariant metric tensor $\boldsymbol{g}(k_i)$. It can be shown that $g=J^2$, where $J$ is the Jacobian (determinant of the Jacobian matrix), and g is the determinant of the co-variant metric tensor. Thus, this functional, similar to the Winslow functional \eqref{winslow}, has a barrier. The value of the functional approaches inf{i}nity when the cells degenerate. That is $\vert\boldsymbol{J}\vert\,\rightarrow\,0$. Thus, this functional produces unfolded grids. Numerical experiments also prove this feature of the functional. The above functional can remove the folded grids produced by the Liao functional. The Modified Liao functional can also be written as follows
\begin{equation}
\mathcal{F}(\mathbf{\mathbf{x},\mathbf{y}}) = \sum_{k=1}^n\left[{\sum_{i=1}^4{(\mathsf{K}(\boldsymbol{J}(k_i)))^2}}\right]\enspace{.}
\label{modlina_new}
\end{equation}
The Modified Liao functional minimizes the square of the condition number where as the Winslow functional minimizes the condition number.
\section{Numerical Examples}
\label{num_examp}
We are interested in f{i}nding such a mesh for which the gradient of the functionals vanish. The minimization of functionals can be performed by a line search algorithms such as Newton's iteration. For the numerical experiments instead of performing the global optimization we solved the local minimization problems for a single node at a time \citep{sanjay4}. In all numerical examples, initial grid was generated by linear transfinite interpolation \citep{gridbook}.
\subsection{Adaptive Grid by Area Functional}
\label{ex_area}
It is generally not recommended to uniformly refine the whole mesh in the hope of capturing the underlying physics. It is desired to adapt a given grid to the requirement of the underlying problem. A grid generating algorithm should be able to allocate more grid nodes in the part of the domain where large solution gradients occur, and fewer grid nodes in the part of the domain where solution is flat. Such grids are called solution-adaptive. Behaviour of the underlying solution can be obtained by posteriori indicators \citep{posteriori}. These indicators can be computed on a coarse mesh, and they can be used to formulate adaptive function $s(x,y)$ in the equation \eqref{area_fun_1}. In the present work, the adaptive function $s(x,y)$ is given in the analytic form. 

\begin{center}
\begin{figure}
\begin{minipage}[t]{6.3cm}
\includegraphics[width=1.3\textwidth]{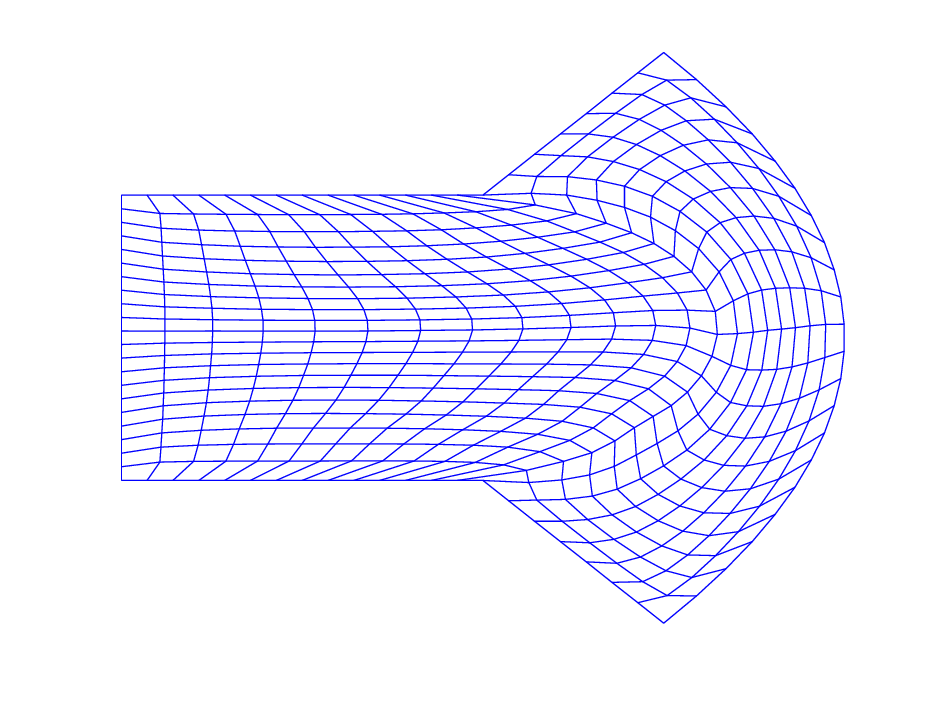}
\caption{Example (\ref{ex_area}) : Adapted Grid by Area Functional. Adapted Functional is give as $s(x,y) = 5.0+200.0\,\vert \sin(2\,\pi\,x)\,\sin(2\,\pi\,y)\vert$.}
\label{fig:Example_surf_dist}
\end{minipage} 
\hfill
\begin{minipage}[t]{6.3cm}
\includegraphics[width=1.3\textwidth]{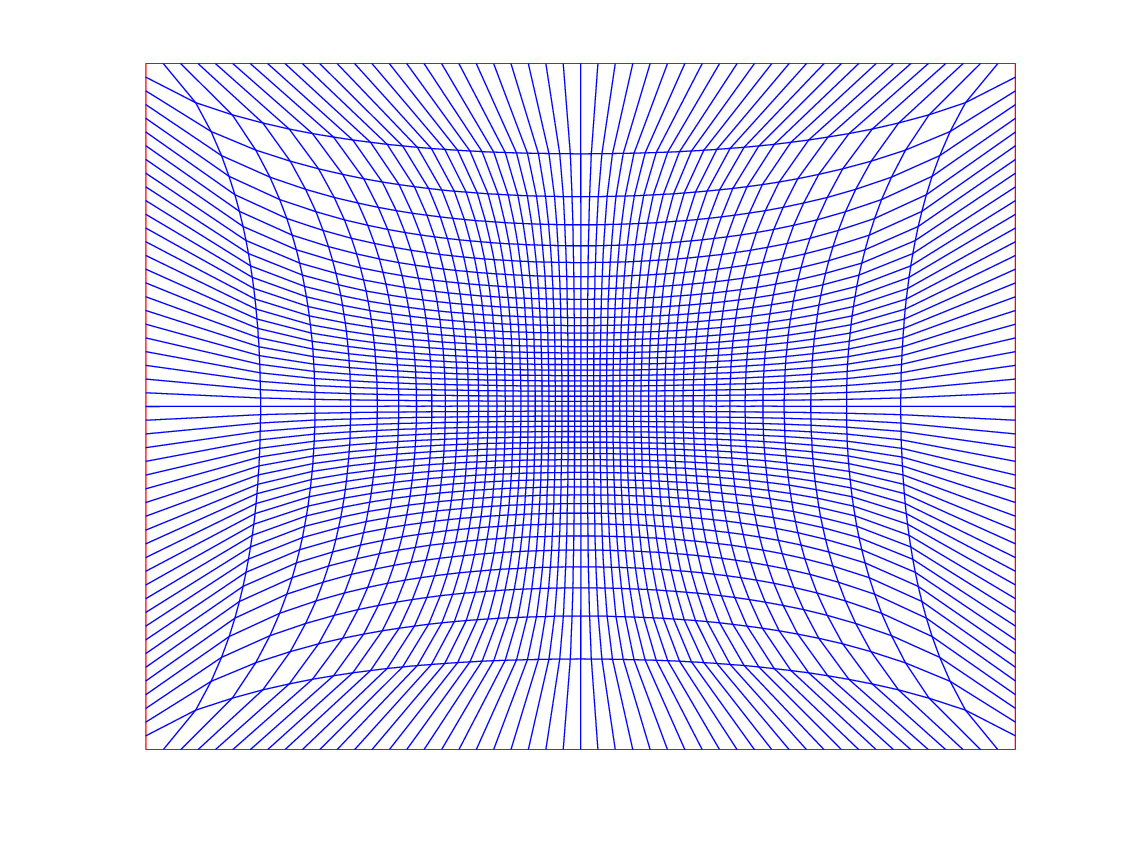}
\caption{Example (\ref{ex_area}) : Adapted Grid by Area Functional. Adapted Functional is give as $s(x,y) = 5.0+200.0\,\left[\sin(2\,\pi\,x)\,\sin(2\,\pi\,y)\right]$.}
\label{fig:Example_velo_1}
\end{minipage}
\hfill
\end{figure}
\end{center}

Figures \ref{fig:Example_surf_dist} and \ref{fig:Example_velo_1} report the outcome of our numerical experiments. It should be noted that even after adaptation the quadrilateral meshes are convex. One other advantage of mesh adaptation by Area functional is that it preserves the mesh topology, and writing a solver for a structured mesh is easier compared to unstructured mesh.
\subsection{Winslow Functional vs Algebraic Method}
\label{examp_2}
\begin{figure}
 \begin{minipage}[t]{8cm}
 \begin{center}
\includegraphics[width=8cm,clip]{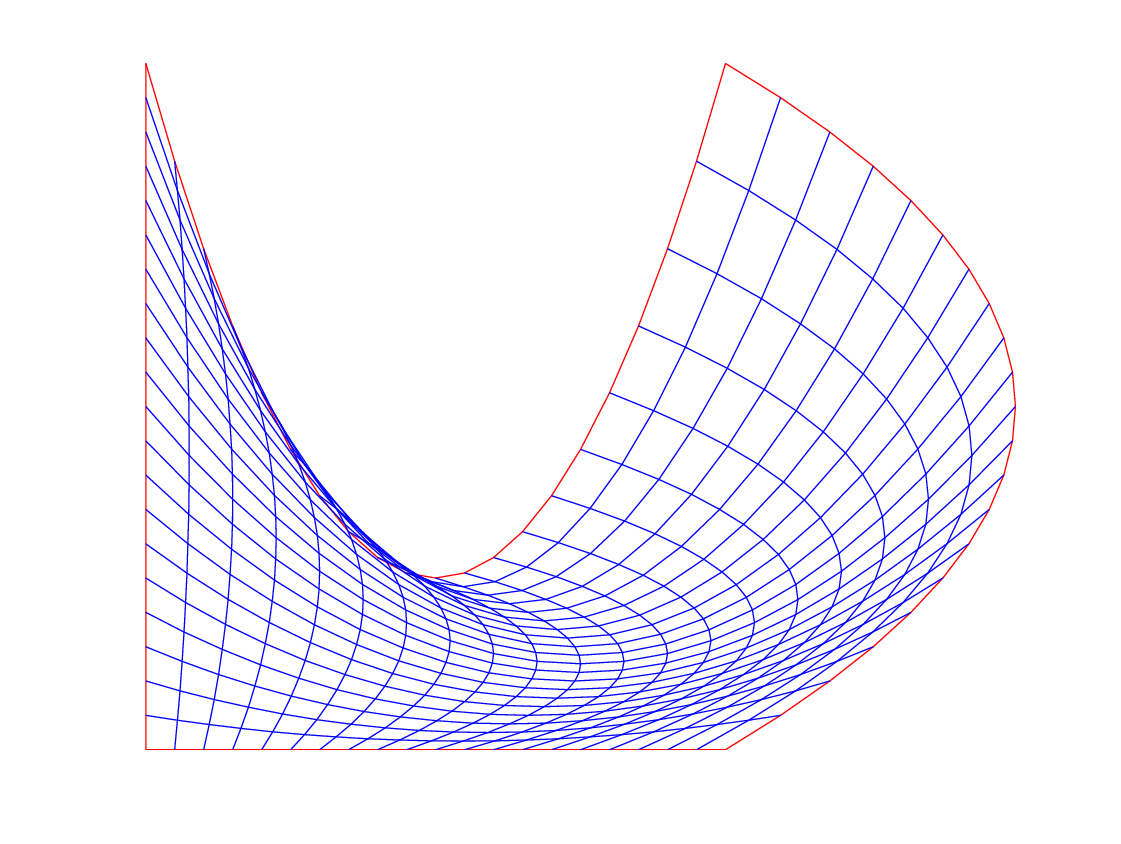}
\caption[Transfinite Interpolation]{\label{labelFig1} Example (\ref{examp_2}) : Folded Grid by Transfinite Interpolation.}
 \end{center}
 \end{minipage}
 \hfill
 \hspace*{0.2cm}
 \begin{minipage}[t]{8cm}
 \begin{center}
 \includegraphics[width=8.0cm,clip]{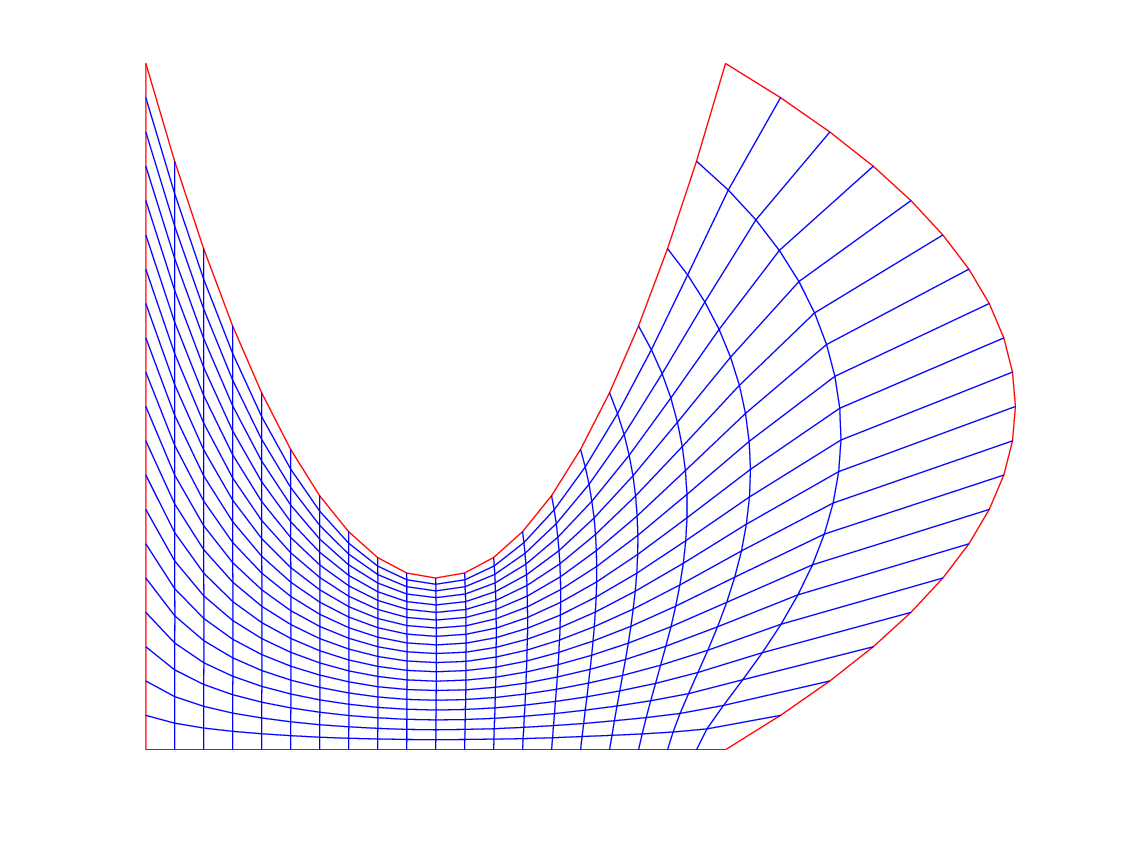}
 \caption[Elliptic Grid]{\label{labelFig2} Example (\ref{examp_2}) : Smooth Grid by Winslow Functional.}
 \end{center}
 \end{minipage}
 \hfill
 \end{figure}

Algebraic grid generation methods such as transfinite interpolations \citep{gridbook} are extensively used for generating grids. Though, they are one of most simplest method of grid generation but algebraic methods can produce folded grids for curved domains as can be seen in the Figure \ref{labelFig1}. One other disadvantage of algebraic grid generation is that boundary discontinuity can prorogate inside the domain. It is clear from Figure \ref{labelFig2} that Winslow functional smooth the grid, and removes the folded grid lines.
\subsection{Liao, Modified Liao and Area Functionals}
\label{example_3}
\begin{figure}
 \begin{minipage}[t]{5cm}
 \begin{center}
\includegraphics[width=5cm,clip]{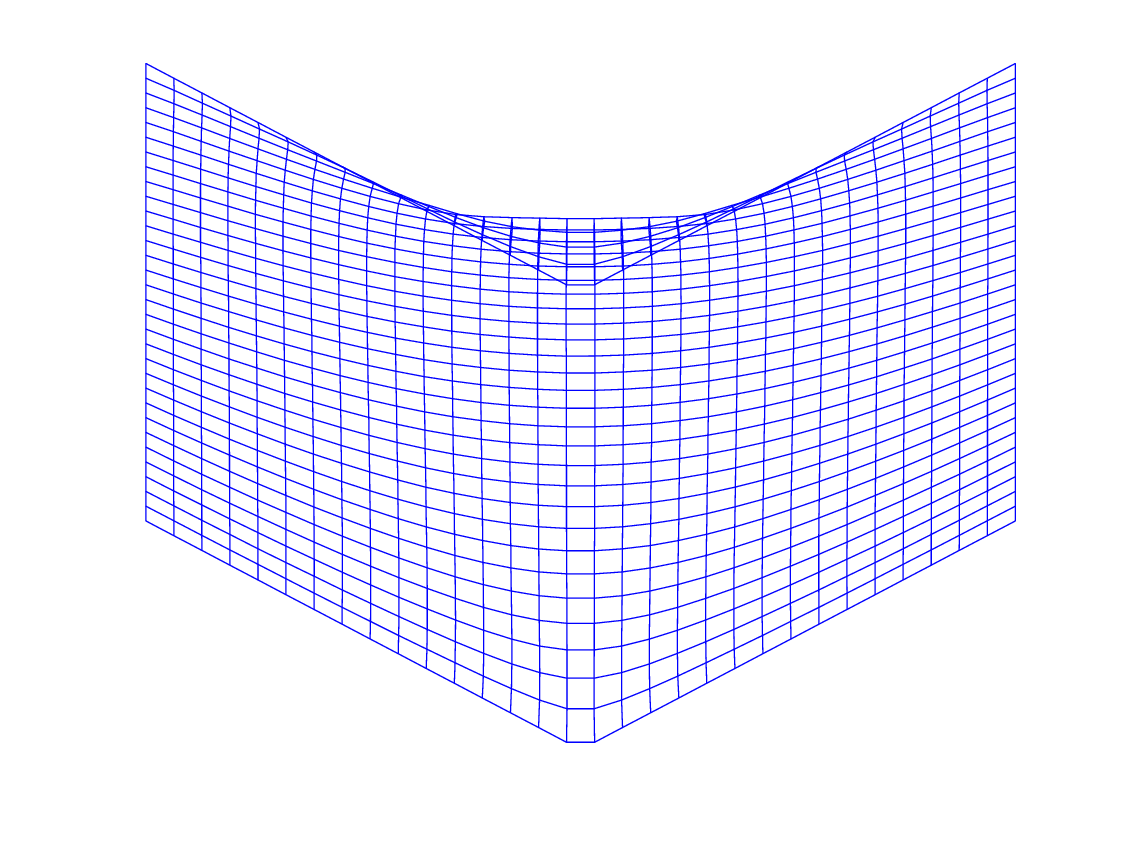}
\caption[Liao Functional]{\label{liao_} Example (\ref{example_3})$\colon$ Folded Grid by the Liao Functional.}
 \end{center}
 \end{minipage}
\hfill
 \begin{minipage}[t]{5cm}
 \begin{center}
 \includegraphics[width=5.0cm,clip]{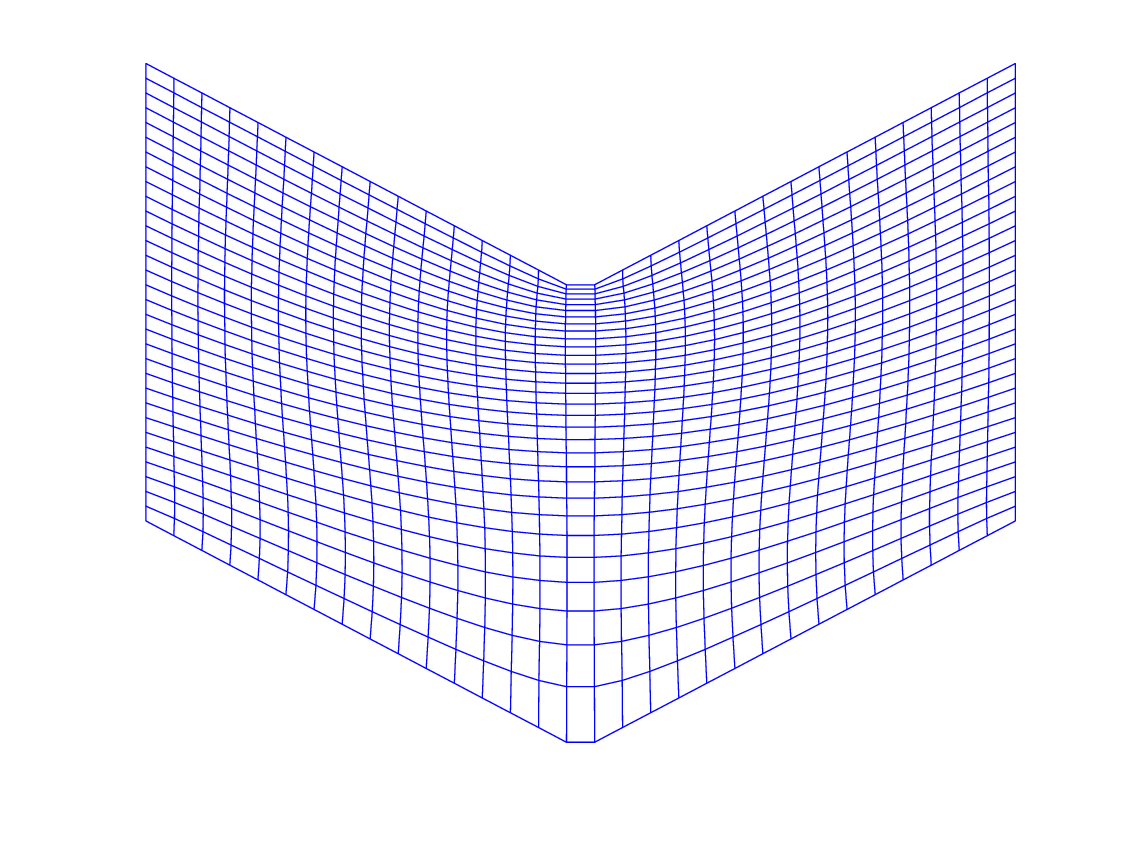}
\caption[Modified Liao Functional]{\label{modliao_} Example (\ref{example_3})$\colon$ Grid by the Modif\mbox{}ied Liao Functional.}
 \end{center}
 \end{minipage}
\hfill
 \begin{minipage}[t]{5cm}
 \begin{center}
 \includegraphics[width=5.0cm,clip]{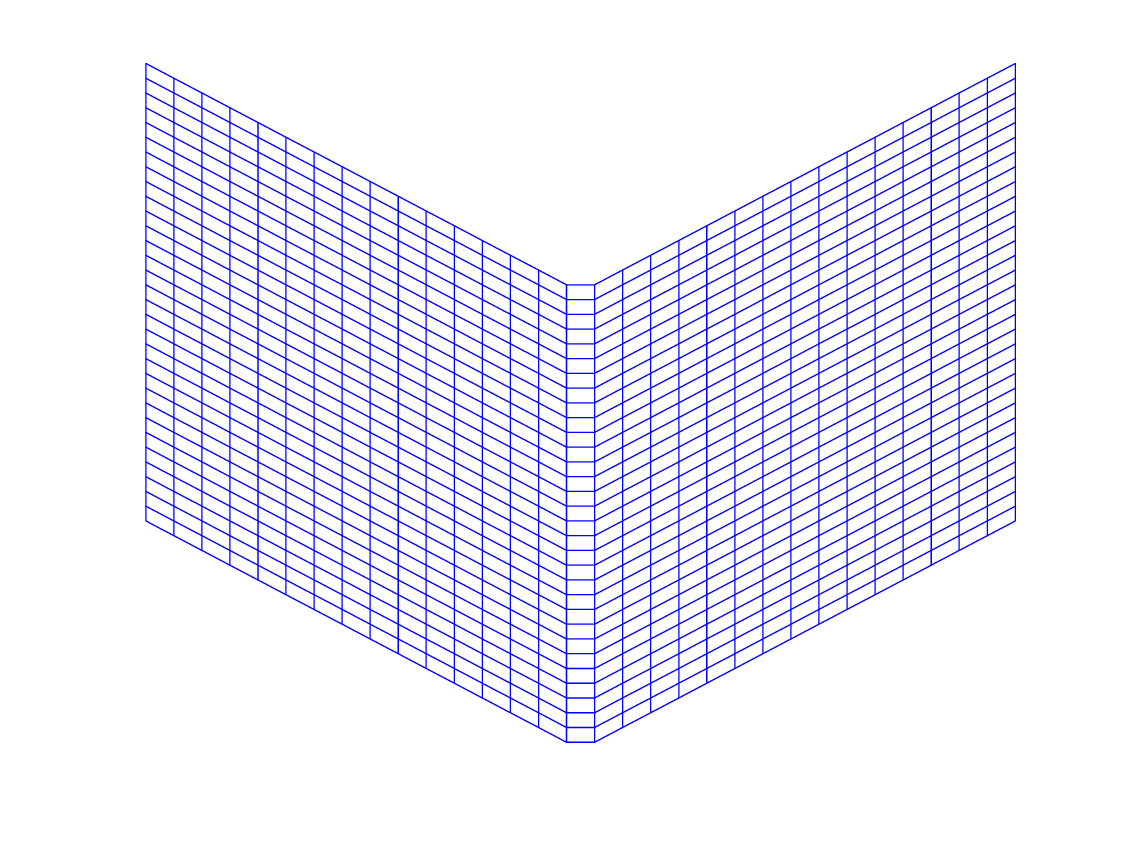}
 \caption[Area Functional]{\label{area_001} Example (\ref{example_3})$\colon$ Grid by the Area Functional with $s(x,y)=1.0$.}
 \end{center}
 \end{minipage}
 \end{figure}

In this example, we perform experiments for comparing Liao, Modified and Area functional on a simple domain. Outcome of our results are shown in Figures \ref{liao_}, \ref{modliao_} and \ref{area_001}. It can be seen from these figures that Modified Liao functional does indeed removes the inverted elements from the mesh but still the quality of the mesh generated by the area functional shown in the Figure \ref{area_001} is certainly better than both Liao and Modified Liao.
\subsection{Length, Area and Knupp's Functionals}
\label{example_4}
\begin{figure}
 \begin{minipage}[t]{5cm}
 \begin{center}
\includegraphics[width=5cm,clip]{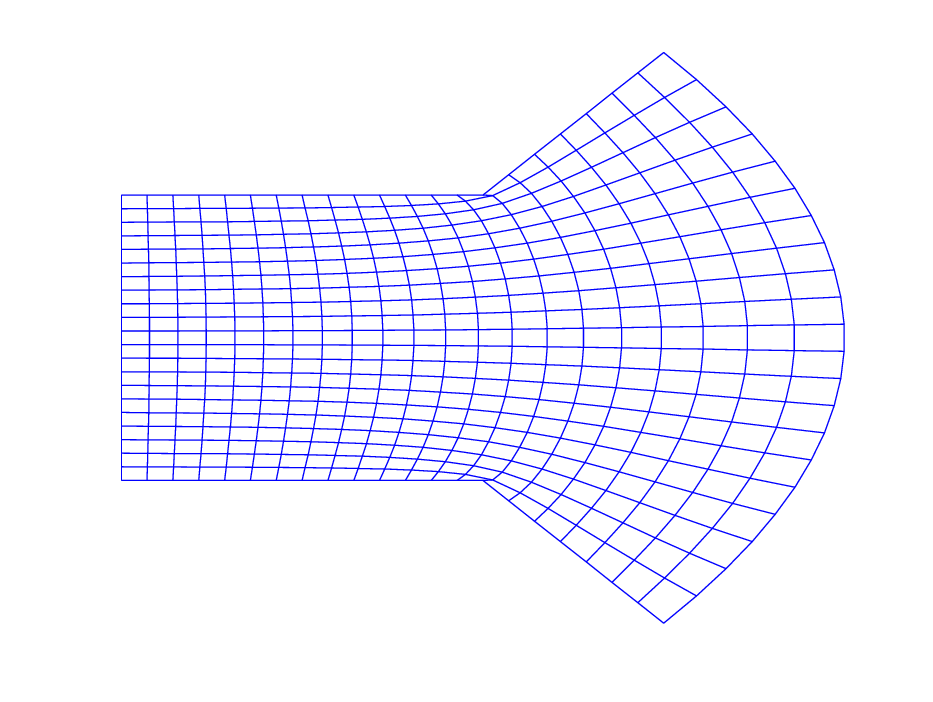}
\caption[Length Functional]{\label{length} Example (\ref{example_4})$\colon$ Grid by the Length functional.}
 \end{center}
 \end{minipage}
 \hfill
 \begin{minipage}[t]{5cm}
 \begin{center}
 \includegraphics[width=5.0cm,clip]{area_00.eps}
\caption[Grid by the Area Functional]{\label{area} Example (\ref{example_4})$\colon$ Grid by the Area Functional with $s(x,y) = 1.0$.}
 \end{center}
 \end{minipage}
 \hfill
 \begin{minipage}[t]{5cm}
 \begin{center}
\includegraphics[width=5.0cm,clip]{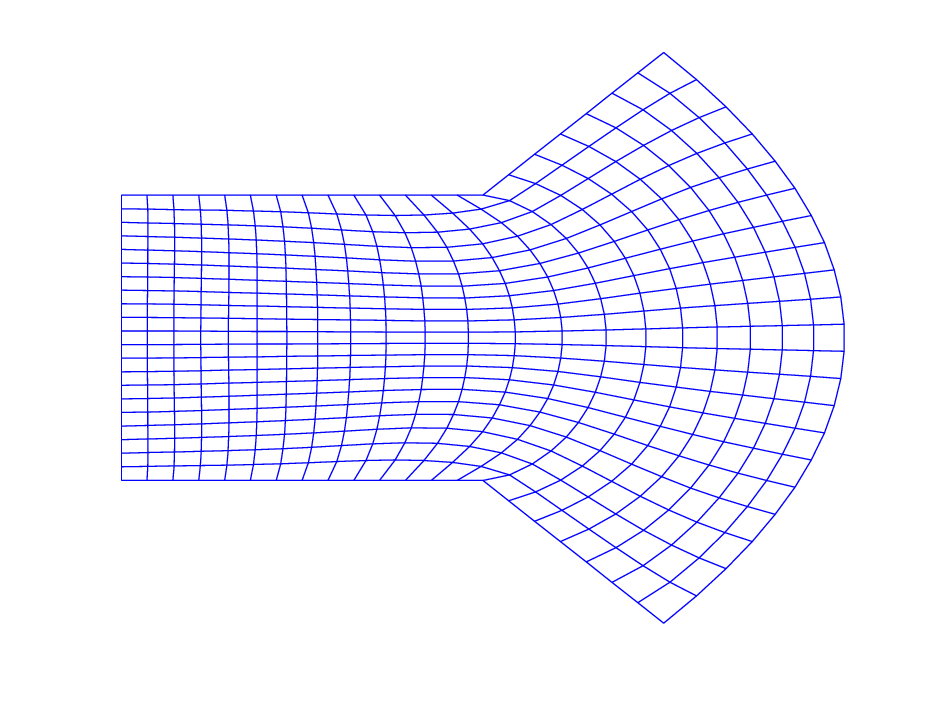}
\caption[Grid by Area and Orthogonality Functional]{\label{area_ortho} Example (\ref{example_4})$\colon$ Grid by the Knupp's functional.}
 \end{center}
 \end{minipage}
 \end{figure}

In this example, we compare the Length, the Area and the Knupp functionals. Figures \ref{length}, \ref{area} and \ref{area_ortho} are the outcome of our numerical work. The Figure \ref{length} is a grid by the Length functional, the Figure \ref{area} is a grid by the Area functional, and the Figure \ref{area_ortho} is a grid by the Knupp's functional. It can be seen that grid by the Area and Knupp's functional are better than the grid produced by the Length functional.
 \subsection{Length and Knupp's Functionals}
\label{example_5}
\begin{figure}
 \begin{minipage}[t]{8cm}
 \begin{center}
\includegraphics[width=8cm,clip]{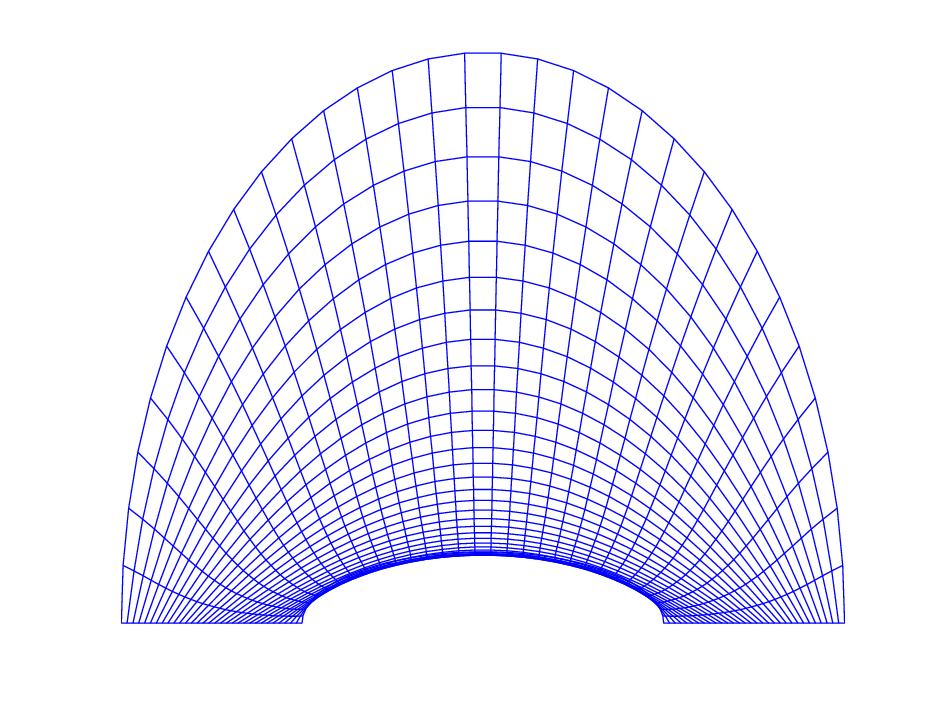}
\caption[Length Functional]{\label{length_2} Example (\ref{example_5})$\colon$ Grid by the Length Functional.}
 \end{center}
 \end{minipage}
 \hfill
 \begin{minipage}[t]{8cm}
 \begin{center}
 \includegraphics[width=8.0cm,clip]{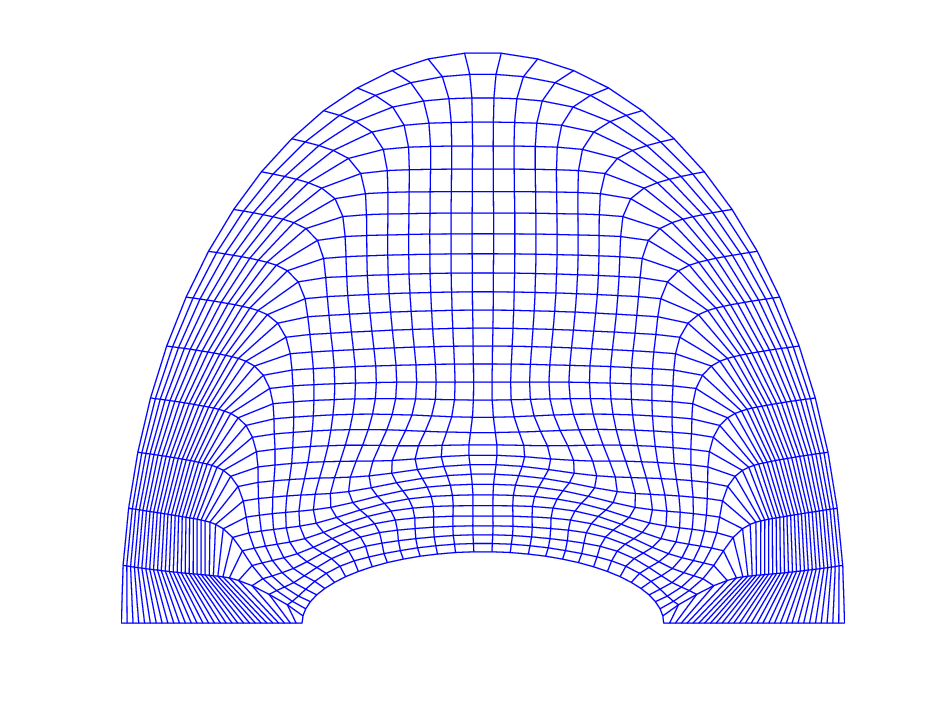}
\caption[Ortho Area]{\label{ortho_area} Example (\ref{example_5})$\colon$ Grid by the Knupp's Functional.}
 \end{center}
 \end{minipage}
 \hfill
 \end{figure}

In this example, we are comparing Length, and the Knupp's Functional. Outcome of our numerical work is reported in Figures \ref{length_2} and \ref{ortho_area}. Figure \ref{length_2} is a grid generated by the Length functional. Figure \ref{ortho_area} is grid generated by the Knupp's functional. It can be seen that the grid generated by the Knupp's functional is of superior quality.
\section{Conclusions}
\label{conclusion_}
We have presented the formulation of various functionals for generating quadrilateral meshes, and an analysis of Winslow and Modified Liao functionals that is consistent with the numerical experiments. Numerical experiments show that Winslow and Modified Liao functionals can remove the folded grids as was expected from theoretical analysis. It has been shown that Area functionals can be used for generating robust adaptive meshes. Further research is required in formulating adaptive function from a posteriori error estimators.

\end{spacing}
\end{document}